\documentstyle [14pt]{article}
\begin{document}
\title{
Free coherent states
\\
and distributions on $p$-adic numbers
}
\author{S.V.Kozyrev}
\maketitle

\footnotetext{
e-mail: kozyrev@class.mian.su,\quad Address:  Inst.Chem.Phys., OSV, 117334,
Kosygina 4, Moscow}

\begin{abstract}
Free coherent states for a system with two degrees of freedom is defined.
A linear map  of the space of free coherent states to the space of
distributions on 2-adic disc  is constructed.
\end{abstract}

\section{Construction of free (or Boltzmannian) coherent states}

Free (or Boltzmannian)
Fock space has been considered in some recent works
on quantum chromodynamics
\cite{MasterFld}, \cite{GopGross}, \cite{DougLi}
and noncommutative probability
\cite{AccLu}, \cite{Maassen}.

The subject of this work is free coherent states.
We will introduce free coherent states and investigate
the  space of coherent states corresponding to a fixed eigenvalue
of the operator of annihilation.
In the paper \cite{qalg}  the  subset  of
the set of free coherent states was constructed.
The main result of the  paper \cite{qalg} is the construction of
the homeomorphism from the ring of integer 2-adic numbers to
the subset of the space of free coherent states with topology defined
by the Hilbert metric.

The result of the present paper is the generalization of the construction
of  the  paper \cite{qalg}.
We introduce the space of coherent states.
The main result of the present paper is the construction of
the linear map of the space of coherent states to
the space of distributions on the ring of integer 2-adic numbers.
Coherent states introduced in  the  paper \cite{qalg} correspond
under this map to $\delta$-functions.

We will consider the system with two degrees of freedom.
The system with one degree of freedom was investigated in \cite{KozAAV}.

The free commutation relations are particular case of
$q$-deformed relations
$$A_i A_j^{\dag}-qA_j^{\dag}A_i=\delta_{ij}  $$
with $q=0$.
A correspondence of  $q$-deformed commutation relations and
non-archimedean (ultrametric) geometry was discussed in
\cite{ArVqGPnAG}.
Non-archimedean mathematical physics
was studied in \cite{VVZ}.

Free coherent states lies in
the free Fock space. Free  (or Boltzmannian) Fock space
$F$ over a Hilbert space  $H$ is the completion of the tensor algebra
$$F=\oplus_{n=0}^{\infty}H^{\otimes n}.$$
Creation and annihilation operators are defined in the following way:
$$
A^{\dag}(f) f_{1}\otimes...\otimes f_{n}=f\otimes
f_{1}\otimes...\otimes f_{n}
$$
$$
A(f) f_{1}\otimes...\otimes
f_{n}=<f,f_{1}> f_{2}\otimes...\otimes f_{n}
$$
where  $<f,g>$ is the scalar product in the Hilbert space $H$.
Scalar product  in the free Fock space is defined by the standard
construction of the direct sum of tensor products of Euclidean spaces.

We consider the case   $H=C\oplus C$, where $C$ is the field of
complex numbers.
In this case we have two creation operators
$A^{\dag}_{0}$, $A^{\dag}_{1}$
and two annihilation operators
$A_{0}$, $A_{1}$ with commutation relations
\begin{equation}\label{aac}
A_{i}A_{j}^{\dag}=\delta_{ij}.
\end{equation}
The vacuum vector $\Omega$  in the free Fock space satisfies
\begin{equation}\label{vacuum} A_{i}\Omega=0.
\end{equation}

We define
free coherent states   in the following way.

Let us consider an infinite sequence of complex numbers
$U=u_{0}u_{1}u_{2}u_{3}...$, $u_i$ are complex numbers.

We introduce the free coherent state $X_{U}$  as
the formal series
$$X_{U}=\sum_{k=0}^{\infty}\lambda^k X_{U}^{k}.$$
Here $X_{0}=\Omega$ is vacuum and
\begin{equation}\label{def_coherent}
X_{U}^{k+1}=\left(u_{k}A_0^{\dag}+(1-u_{k})A_1^{\dag}\right)X_{U}^{k}.
\end{equation}
This formal series define a functional with dense domain in the free Fock space.

Free coherent states are formal eigenvectors of the annihilation operator
$$
A_{0}+A_{1};
$$
for the eigenvalue $\lambda$, i.e.
$$
(A_{0}+A_{1})X_{U}=\lambda X_{U}\qquad \forall U.
$$
We define
the linear space  $X'$ of free coherent states (formal eigenvectors of
$A_{0}+A_{1}$) as the linear span of functionals $X_U$.

\section{Construction of the isometry of a subset of the space
of coherent states to 2-adic disc}

In the present section the eigenvalue $\lambda\in (0,1)$.
The subject of the present section is the investigation of
the set of coherent states $X_U$ where $U$ runs on all
infinite sequences of 0 and 1.
Degeneration of this eigenspace is parametrized by
the set $\{U\}$ of infinite sequences of 0 and 1.

These sequences have  a natural interpretation as 2-adic numbers.
Every sequence $U$ is in one-to-one correspondence with a 2-adic number
$U=\sum_{i=0}^{\infty}u_i 2^{i}$.

Let us consider the  metric $\rho$ of free Fock space on the set of free
coherent states.
Let $U$, $V$ be arbitrary sequences of 0 and 1. These sequences
coincide up to an element with number $k-1$. 2-adic distance
of corresponding 2-adic numbers is $||U-V||_2=2^{-k}$.
Let $X_{U}$, $X_{V}$ be corresponding free coherent states.
We have:
$$
\rho(X_{U},X_{V})^2=2\sum_{i=k}^{\infty}\lambda^{2i}=
\lambda^{2k}\frac{2}{1-\lambda^{2}}.
$$
Therefore the metric  $\rho$ is proportional to the degree of 2-adic metric.

The metric  $\rho$ is an ultrametric.
This means that the metric  $\rho$ obeys the strong triangle inequality
$$
\rho(X_{U},X_{V})\le \max(\rho(X_{U},X_{W}),\rho(X_{V},X_{W}))
$$
for arbitrary $U$, $V$, $W$.

We get the following proposition.

{\bf   Proposition.}{\sl\qquad
The map $X_U\mapsto \sum_{i=0}^{\infty}u_{i}2^{i}$
is the
isometry  of the set of free coherent states  $\{ X_{U} \}$  parametrized by
infinite sequences $U$ of 0 and 1
with the metric $\rho$ of scalar product
to  the ring of integer 2-adic numbers with the metric proportional
to the degree of $||U,V||_2$.
}

\section{Construction of the map of  the space
of coherent states to the space of distributions on the 2-adic disc}

In the present section the linear map of the space $X'$
of coherent states to the space of distributions on the 2-adic disc
will be constructed. In the present section we take $\lambda\in (0,\sqrt{2})$.

Distributions on the 2-adic disc are linear functionals on the space
of locally constant functions \cite{VVZ}.
Therefore we have to construct the coherent state that corresponds to a locally
constant function.

Let us consider the sequence $W=w_0 w_1 w_2...$, $w_i$ are arbitrary
complex numbers.
Let us introduce the sequence
\begin{equation}\label{trun}
W_k=w_0 w_1 w_2...w_{k-1}\frac{1}{2}\frac{1}{2}\frac{1}{2}...
\end{equation}
For $\lambda\in (0,\sqrt{2})$ the coherent state $X_{W_k}$ lies in the
Hilbert space (the correspondent functional is bounded).

We will prove that
an arbitrary locally constant function  corresponds to the
linear combination of coherent
states $X_{W_k}$.
The linear span of
the coherent states $X_{W_k}$
we will denote by $X$.
Every vector in  $X$ is a function of $\lambda$. We will investigate
the  properties of the space $X$ with the scalar product
\begin{equation}\label{product_trun}
<X_{U_i},X_{V_j}>=
\lim_{\lambda\to\sqrt{2}-0}\left(1-\frac{\lambda^2}{2}\right)
(X_{U_i},X_{V_j}).
\end{equation}

It is sufficient to find coherent states that correspond to locally constant
functions of the type
$$
\theta_k(x-x_0)=\theta(2^{k}||x-x_0||_2);\quad
\theta(t)=0, t>1;\quad \theta(t)=1, t\le 1.
$$
Here $x$, $x_0\in Z_2$ lies in
the ring of integer 2-adic numbers and the function $\theta_k(x-x_0)$ equals to
1 on the disc $D(x_0,2^{-k})$ of radius $2^{-k}$ with the center in $x_0$ and
equals to 0 outside this disc.

Let us consider the sequences $U=u_0 u_1 u_2...$, $u_i=0,1$
and $V=v_0 v_1 v_2...$, $v_i=0,1$
that correspond to 2-adic numbers $U=\sum_{i=0}^{\infty}u_i 2^i$
and $V=\sum_{i=0}^{\infty}v_i 2^i$.
We have the following lemma.

{\bf Lemma 1.}\qquad
{\sl
The  limit of
the scalar product of coherent states $X_{U_k}$, $X_{V_l}$
equals to the integral on 2-adic disc with respect to the Haar measure
$$
\lim_{\lambda\to\sqrt{2}-0}\left(1-\frac{\lambda^2}{2}\right)
(X_{U_k},X_{V_l})=
$$
$$
=\frac{1}{\mu(D(U,2^{-k}))\mu(D(V,2^{-l}))}
\int_{Z_2}\theta_k(x-U)\theta_l(x-V)dx.
$$
Here $\mu(D)$ is the Haar measure of the disc $D$.

}

{\bf Proof}

The formula (\ref{def_coherent}) for the coherent states
$X_{U_k}$ and $X_{V_l}$ have the following form
$$
X_{U_k}=\sum_{i=0}^{k-1}\lambda^i X_{U_k}^{i}+
\lambda^k\sum_{i=0}^{\infty}\lambda^i
\left(\frac{1}{2}A_0^{\dag}+\frac{1}{2}A_1^{\dag}\right)^i X_{U_k}^{k};
$$
$$
X_{V_l}=\sum_{i=0}^{l-1}\lambda^i X_{V_l}^{i}+
\lambda^l\sum_{i=0}^{\infty}\lambda^i
\left(\frac{1}{2}A_0^{\dag}+\frac{1}{2}A_1^{\dag}\right)^i X_{V_l}^{l}.
$$
Let $k\le l$.

If the first $k$ indices of the sequence $U$  coincide with
the first $k$ indices of the sequence $V$ then the scalar product
have a following form
$$
(X_{U_k},X_{V_l})=\sum_{i=0}^{k-1}\lambda^{2i}+
\sum_{i=k}^{l-1}\lambda^{2i}\left(\frac{1}{2}\right)^{i-k}+
\sum_{i=l}^{\infty}\lambda^{2i}\left(\frac{1}{2}\right)^{i-k}.
$$
If the first $k$ indices of the sequence $U$  do not coincide with
the first $k$ indices of the sequence $V$ then the series for scalar product
$(X_{U_k},X_{V_l})$ contains only the finite number of terms.
Therefore the limit
$\lim_{\lambda\to\sqrt{2}-0}\left(1-\frac{\lambda^2}{2}\right)(X_{U_k},X_{V_l})$
equals to $\min(2^{k},2^{l})$ if one of the discs $D(U,2^{-k})$
and $D(V,2^{-l})$  contains another and equals to 0 if these discs
do not intersect.

Therefore the  limit of the scalar product of vectors $X_{U_k}$,
$X_{V_l}$ in the free Fock space corresponds to the integral on 2-adic disc
with respect to the Haar measure. Coherent state $X_{U_k}$ corresponds to
the locally constant function   $\frac{1}{\mu(D(U,2^{-k}))}\theta_k(x-U)$.

Let us investigate functionals on the space $X$ of locally constant functions.

Coherent states $X_U$ where
sequences $U$ consist of 0 and 1  correspond to  $\delta$-functions.

Let $U$ and $V$ be sequences of 0 and 1.
We have the following lemma.

{\bf Lemma 2.}\qquad
{\sl
The limit of  the action of the functional $X_{U}$
on the vector  $X_{V_j}$
have the following form
$$
\lim_{\lambda\to\sqrt{2}-0}\left(1-\frac{\lambda^2}{2}\right)
(X_{U},X_{V_j})=
\frac{1}{\mu(D(V,2^{-j}))}\int_{Z_2}\delta(x-U)\theta_j(x-V)dx.
$$
}

Therefore the coherent state $X_{U}$ corresponds to
the $\delta$-function   $\delta(x-U)$.

{\bf Lemma 3.}\qquad
{\sl
Vectors $X_{W_k}$ lie in the domain of the functional $X_{U}$
for arbitrary sequence $U$ of complex numbers.
}

{\bf Proof}

The formula (\ref{def_coherent}) for the coherent state
$X_{W_k}$ have the following form
$$
X_{W_k}=\sum_{i=0}^{k-1}\lambda^i X_{W_k}^{i}+
\lambda^k\sum_{i=0}^{\infty}\lambda^i
\left(\frac{1}{2}A_0^{\dag}+\frac{1}{2}A_1^{\dag}\right)^i X_{W_k}^{k}.
$$
The action of the functional $X_{U}$ is defined by the following formal series
\begin{equation}\label{action}
(X_{U},X_{W_k})=\sum_{i=0}^{\infty}\lambda^{2i} (X_{U}^{i},X_{W_k}^{i}).
\end{equation}
Let us calculate
$$
(X_{U}^{k+i},X_{W_k}^{k+i})=(X_{U}^{k+i-1},
\left(u^{*}_{k+i-1}A_0+(1-u^{*}_{k+i-1})A_1\right)
\left(\frac{1}{2}A_0^{\dag}+\frac{1}{2}A_1^{\dag}\right)
X_{W_k}^{k+i-1})
$$
for $i>0$. We have
$$
\left(u^{*}_{k+i-1}A_0+(1-u^{*}_{k+i-1})A_1\right)
\left(\frac{1}{2}A_0^{\dag}+\frac{1}{2}A_1^{\dag}\right)=
\frac{1}{2}(u^{*}_{k+i-1}+1-u^{*}_{k+i-1})=\frac{1}{2}.
$$
Therefore
$$
(X_{U}^{k+i},X_{W_k}^{k+i})=\frac{1}{2}(X_{U}^{k+i-1},X_{W_k}^{k+i-1});
$$
we have $\frac{\lambda^2}{2}<1$
and the series (\ref{action}) converges.

Therefore an arbitrary coherent state $X_{U}$ corresponds to a distribution
on the 2-adic disc.

As a corollary of the Lemma 3 we introduce the linear map $\phi$
$$
\phi:\quad X'\to D'(Z_2);
$$
$$
\phi:\quad X\to D(Z_2);
$$
where $D'(Z_2)$ is the space of distributions on the 2-adic disc $Z_2$
and   $D(Z_2)$ is the space of locally constant functions on the 2-adic disc
$Z_2$. We have the following corollary of Lemma 3.

{\bf Corollary.}\qquad
{\sl
The formula
$$
\lim_{\lambda\to\sqrt{2}-0}\left(1-\frac{\lambda^2}{2}\right)
(X_{U},X_{W_k})
$$
defines the action of the distribution $\phi(X_{U})$ on the locally constant
function $\phi(X_{W_k})$.
}

We get the following theorem.

{\bf Theorem.}\qquad
{\sl
The map $\phi$
$$
X_{V_j}\mapsto \frac{1}{\mu(D(V,2^{-j}))}\theta_j(x-V)
$$
where $V$ is an arbitrary sequence of 0 and 1, $j=0,1,2...$
extends to the isomorphism of the space $X$  of coherent states
of a type (\ref{trun})
to the space $D(Z_2)$ of locally constant functions on the ring of integer
2-adic numbers with the scalar product defined by the Haar measure.

The scalar product of locally constant functions  in $L_2$
with respect to the Haar measure
equals to the limit of the scalar product
in the free Fock space
of corresponding
coherent states
$$
\lim_{\lambda\to\sqrt{2}-0}\left(1-\frac{\lambda^2}{2}\right)
(X_{U_i},X_{V_j})=
$$
$$
=\frac{1}{\mu(D(U,2^{-i}))\mu(D(V,2^{-j}))}
\int_{Z_2}\theta_i(x-U)\theta_j(x-V)dx;
$$
here $U$, $V$ are arbitrary sequences of 0 and 1.

An arbitrary coherent state $X_{U}$
corresponds to a distribution on the ring of integer 2-adic numbers
with the action on locally constant functions defined by the formula
$$
(\phi(X_{U}),\phi(X_{W_k}))=
\lim_{\lambda\to\sqrt{2}-0}\left(1-\frac{\lambda^2}{2}\right)
(X_{U},X_{W_k}).
$$

}
%
%
%
%
%

\section{Conclusion}

Coherent states introduced in the present paper can be interpreted
in the spirit of noncommutative geometry.

For example let us consider the real plane $R^2$ with coordinates $x$, $y$.
Coordinates $x$, $y$ we will consider as operators acting on functions on $R^2$.
Let us introduce ''coherent states'' $f$ as eigenvectors of the operator
$x+y$:
$$
(x+y)f=f.
$$
It is easy to see that these eigenvectors are distributions with support
on the line $x+y-1=0$.

An arbitrary algebraic manifold  with equation $F(x,y)=0$
can be described in the same way. The distribution $\eta(x,y)$ on this manifold
is a distribution on the plane $R^2$ that is annihilated by the polynom
$F(x,y)$:
$$F(x,y)\eta(x,y)=0.$$
A point on the algebraic manifold with equation $F(x,y)=0$
will correspond to the $\delta$-function with support in the point
$(x_0,y_0)$ such that $F(x_0,y_0)=0$.

For the case of free (Boltzmannian) coherent states we have
a noncommutative algebraic manifold with equation $A_0+A_1-\lambda=0$.
The result of the present paper can be interpreted in the following form:
the limit for $\lambda\to \sqrt{2}-0$
of the manifold with equation $A_0+A_1-\lambda=0$ is 2-adic disc.

This result illustrates the conjecture of the paper \cite{ArVqGPnAG}:
noncommutative algebrae can be used for a deformation of real manifolds
to $p$-adic manifolds.

\vspace{3mm}
{\bf Acknowledgments}

Author is grateful to I.V.Volovich for discussions.

\end{document}